\documentclass[article]{amsart}

\usepackage{verbatim}
\usepackage{amssymb}
\usepackage{amsbsy}
\usepackage{amscd}
\usepackage{amsmath}
\usepackage{amsthm}
\usepackage[mathscr]{eucal}

\newtheorem{theorem}{Theorem}

\newtheorem{problem}{Problem}
\newtheorem{remark}{Remark}
\newtheorem{defn}{Definition}\numberwithin{defn}{section}

\numberwithin{figure}{section}
\numberwithin{table}{section}
\numberwithin{remark}{section}
\numberwithin{equation}{section}
\numberwithin{theorem}{section}
\numberwithin{Lemma}{section}
\numberwithin{problem}{section}

\usepackage{graphicx}

\usepackage{subcaption}
\def\Xint#1{\mathchoice
	{\XXint\displaystyle\textstyle{#1}}%
	{\XXint\textstyle\scriptstyle{#1}}%
	{\XXint\scriptstyle\scriptscriptstyle{#1}}%
	{\XXint\scriptscriptstyle\scriptscriptstyle{#1}}%
	\!\int}
\def\XXint#1#2#3{{\setbox0=\hbox{$#1{#2#3}{\int}$ }
		\vcenter{\hbox{$#2#3$ }}\kern-.6\wd0}}

\def\dashint{\Xint-}

\title[Biomembranes]{Small deformations of spherical biomembranes}

\author[C.M. Elliott, L. Hatcher and P.J. Herbert]{C.M. Elliott, L. Hatcher and P.J. Herbert}

\address{Mathematics Institute, University of Warwick}

 \email{c.m.elliott@warwick.ac.uk}

\numberwithin{equation}{section}

\subjclass[2010]{}


\begin{document}

\begin{abstract}
In this contribution to the proceedings of the 11th Mathematical Society of Japan (MSJ)  Seasonal Institute (July 2018) we give an overview of  some recent work on a 
mathematical model for small deformations of a spherical membrane. The idea is to consider perturbations to minimisers 
of a surface geometric energy. The model is obtained   from consideration of  second order approximations to a perturbed energy. 
In particular, the considered problems involve particle constraints and surface phase field energies.

\end{abstract}

\maketitle

\section{Introduction}

Mathematical models in cell biology involving biomembranes  are of burgeoning interest within biophysical  applications,  analysis of geometric partial differential equations and numerical analysis.   Biomembranes are thin lipid bilayers which surround cells
and some organelles providing  a barrier between their interior and exterior. They are
principally composed of  lipid molecules whose tails are hydrophobic  and thus forming a bilayer  in which the hydrophilic lipid heads point out of the bilayer. Within the bilayer and its local environment there are  a large number of proteins. 
These proteins have a great variety of forms
and purposes. They may be attached to a surface of the bilayer, embedded into
it or span the entire bilayer causing
local deformations. The great variety of proteins and their possible interactions
with the membrane mean there are many models for membrane deformation. They are  known to regulate
cell morphology and a variety of cellular functions see,
for example, \cite{McMGal05}. 

In
biophysics, established continuum models  treat the biomembrane as a deformable 
fluid surface of infinitesimal thickness whose deformation  is described by an energy functional that depends on the curvature of the surface. 
Our focus is on equilibrium configurations for biomembranes associated with the 
interaction of proteins and the  curvature of the lipid bilayers.  Proteins can induce curvature by
local shaping but also  elastic forces may   organize membrane
proteins. 

In this contribution we consider two models for small perturbations of the surface energy. 
   First we observe that  a large area of study is the 
effect of embedded inclusions which deform
the membrane by their shape \cite{McMGal05, EllGraKor16,  Kie18}. The cross-section of membrane proteins and the
membrane surface usually differs dramatically in scale.  Here we  consider individual proteins  as imposing point constraints or forces on the membrane.
 This leads to fourth order equations with point constraints.
 Secondly we  consider varying concentrations of proteins.  The impact on the membrane morphology is then modelled by concentration-dependent
material properties, phase-dependent spontaneous curvature, and a gradient energy. For this regime we consider continuous distributions of particles 
 which separate the membrane into two phases with differing properties. 


\section{Canham-Helfrich models}
A classical energy for the biomembrane 
is the  energy functional \cite{Can70,Hel73},
\begin{equation}\label{eq:simpleHelfrich}
\mathcal{E}_{B}(\Gamma) := \int_\Gamma \left(\frac{1}{2}\kappa (H-H^s)^2 + \sigma + \kappa_G K \right){\rm d}\Gamma.
\end{equation}
Here $\Gamma$ is a two dimensional hypersurface in $\mathbb R^3$ modelling the biomembrane.  The interior of the biomembrane is $\Omega$, an open bounded 
subset of $\mathbb R^3$, whose boundary $\partial \Omega=\Gamma$ has outward pointing unit normal $\nu$.  We denote by $\mathcal A(\Gamma)$ and $\mathcal V(\Gamma)$ the area and volume of the surface and enclosed bulk domain which satisfy the isoperimetric inequality
$$\mathcal A^3(\Gamma)\ge 36\pi\mathcal V^2(\Gamma).$$
The mean and Gaussian curvatures are denoted by $H$ and $K$. We use the notation $\nabla_\Gamma$ to denote the tangential gradient, 
and $\mathcal H=\nabla_\Gamma\cdot \nu$ to denote the extended Weingarten map with $H$ equal to the sum and $K$ being the product of the two principal curvatures  which are the eigenvalues of $\mathcal H$ associated with two tangential eigenvectors. Note that this definition of the mean curvature differs from the more common definition being the arithmetic mean.
The constants $\kappa,\kappa_G >0$ and $\sigma \geq 0$ measure the bending rigidity and the surface tension of the membrane,
respectively.
The remaining parameter $H^s$ is the spontaneous curvature, that is the preferred curvature of the membrane. This may be due to asymmetry in the configuration of the
lipid bilayer. It may vary over the membrane. Without the surface tension term this energy functional is called the Canham-Helfrich functional. This was introduced by Canham \cite{Can70} in the case $H^s=0$. The spontaneous curvature $H^s$ was introduced by Helfrich \cite{Hel73}. 

The biological context has led to the minimisation problem for these energy functionals being studied with respect to various constraints. Some typical examples would be: taking into account osmotic pressure effects, which leads to a fixed volume constraint; taking the membrane to be incompressible, which leads to an area constraint; supposing that molecules cannot be exchanged between the two layers of the bilayer, which leads to a constraint on the curvature \cite{Sei97}. In the presence of proteins and  varying lipid concentrations these equations are typically coupled with nonlinear partial
differential equations describing the evolution and location of interfaces, concentrations and particles on the surface.

Difficulties in studying these minimisation problems arise since the Euler-Lagrange equations associated
with these surface  energy functionals usually are fourth-order, nonlinear (geometric) partial
differential equations. For example, some minimisation problems of the  Willmore energy were considered in \cite{BauKuw03,Sim93} and some minimisation problems for the Helfrich energy were considered in \cite{ChoMorVen13,ChoVen13,Hel15}.

In  the case that $\kappa_G$ is constant we have by the Gauss-Bonnet theorem
$$\kappa_G\left (\int_\Gamma K +\int_{\partial \Gamma} H_g\right)=\kappa_G 2\pi \chi(\Gamma),$$
where $H_g$ is the geodesic curvature of $\partial \Gamma$ and $\chi(\Gamma)$ is the Euler characteristic. Note that if $\partial \Gamma=\emptyset$ then the corresponding term in the formula vanishes.  If 
either $\kappa_G$ is zero or the integral of the geodesic curvature over $\partial \Gamma$ is constant then we may neglect the Gaussian curvature term in the energy.  The Willmore 
functional
\begin{equation}
\mathcal{W}(\Gamma) :=  \int_{\Gamma} \frac{1}{2}H^2  {\rm d}\Gamma,
\end{equation}
is a special case, well studied in differential geometry, \cite{Wil93}.
 
 Our base functional in this discussion is obtained by setting  $H^s = 0$ and supposing the Gaussian curvature term is constant. The resulting  energy
is
\begin{equation}\label{Willmore+tension}
\mathcal{E}(\Gamma)=\int_{\Gamma}\left(\frac{\kappa}{2}H^2+\sigma\right) {\rm d}\Gamma.
\end{equation}

Let $\Gamma_0$ be a fixed reference compact hypersurface over which we seek a map $X:\Gamma_0\rightarrow \Gamma\subset \mathbb R^3$
that is
$$\Gamma=\{x\in \mathbb R^3~|~x=X(p), p\in \Gamma_0\}.$$ 
The reference hypersurface $\Gamma_0$  might be chosen for convenience of the problem under consideration. In this article we wish to use the  perturbation approach of   \cite{EllFriHob17-a} which  takes advantage of the evaluation of the first variation of the energy functional at $\Gamma_0$ being  zero in order to derive quadratic energy functional approximations.  This allows the study of approximate models which have a simpler variational structure. 
 In the following sections we choose $\Gamma_0$ to be a sphere
$ \mathbb{S}^2(R)$ of radius $R$ with enclosed volume $V(R)$. 

 We  conclude this section by setting out  briefly several problems. In each we constrain the volume enclosed by the surface.
We write, with an abuse of notation. 
$$\mathcal E(X):= \int_{\Gamma_0} \left(\frac{1}{2}\kappa H^2 + \sigma\right){\rm d}\mu$$
where ${\rm d}\mu$ denotes the surface area element of $\Gamma$ associated with the map $X$.
\subsection{Point constraints}
We begin with point constraints which model the interaction of a biomembrane with proteins. To begin we suppose that there is an interaction between particles external to the membrane with a quadratic energy.
\begin{problem}\label{soft}
Let $\{z_j\}_{j=1}^L$ be a collection of disjoint points in $\mathbb R^3$.  
Find $X^*\in K:=\{X\in C^2(\Gamma_0)~|~\mathcal{V}(X(\Gamma_0)) = V_0\}$ such that  

$$\mathcal E(X^*)+\sum_{j=1}^L \frac{1}{2\delta_j}d(X^*(\Gamma_0),z_j)^2=\inf_{X\in K}\left (\mathcal E(X) +\sum_{j=1}^L \frac{1}{2\delta_j}d(X(\Gamma_0),z_j)^2\right ),$$
where $\delta_j >0$ might represent the  reciprocal  of an  elastic  spring constant associated to the particle attachment.
\end{problem}
Here $d\colon \mathbb{R}^3 \times \mathbb{R}^3 \to \mathbb{R}^+$ is the standard metric on $\mathbb{R}^3$, with the usual extension to sets,
\[
	d(X(\Gamma_0),z):= \inf_{y\in X(\Gamma_0)} d(y,z).
\]
In the case of vanishing $\delta_j$ we obtain  point Dirichlet constraints.
\begin{problem}\label{hard}
Let $\{z_j\}_{j=1}^L$ be a collection of disjoint points.  
Find $X^*\in K(\Gamma):=\{X\in C^2(\Gamma_0)~|~\exists p_j\in \Gamma_0 : X(p_j)=z_j, j=1,...,L, ~\mathcal{V}(X(\Gamma_0)) = V_0\}$ such that  

$$\mathcal E(X^*)=\inf_{X\in K}\mathcal E(X).$$
\end{problem}

 The point constraints in Problem \ref{hard} are said to be hard constraints. The energy of Problem \ref{soft} may be viewed as involving   soft  constraints corresponding to the penalisation of  the hard constraints.

\subsection{Two phase surfaces}

 We will consider the case that some intercalated particles within the membrane induce a spontaneous curvature. To model this we introduce an order parameter $\phi:\Gamma\to\mathbb{R}$ to represent the difference of two  relative concentrations, which we assume satisfies the constraint $\dashint_\Gamma \phi =\alpha$ ($\dashint_\Gamma:=\frac{1}{|\Gamma|}\int_\Gamma$). Two phases of the membrane can then be distinguished; a particle rich phase and a particle deficient one. We also introduce a Ginzburg-Landau energy functional to approximate the line tension between the differing phases and obtain the energy functional
\begin{equation}\label{GLhelfrich}
\mathcal E_{JL}^{GL}(X,\phi)=\int_{\Gamma_0}\left (\frac{\kappa}{2}\left(H-H^s(\phi)\right)^2+\sigma+\frac{b\epsilon}{2}|\nabla_\Gamma\phi|^2+\frac{b}{\epsilon}W(\phi) \right){\rm d}\mu.
\end{equation}
Here $W(\phi)$ is a double well potential, for example  given by $W(\phi)=\frac{1}{4}(\phi^2-1)^2$, which favours the two distinct values $\phi=\pm 1$ generating phases, $\epsilon\in\mathbb{R}$ is a smallness parameter commensurate with the width  of a diffuse  interface, $b\in\mathbb{R}$ is a line tension coefficient and as before ${\rm d}\mu$ denotes the surface area element of $\Gamma$ associated with the map $X$. If we take a phase-dependent spontaneous curvature of the form
\begin{displaymath}
H^s(\phi)=\Lambda\phi,
\end{displaymath} 
where $\Lambda\in\mathbb{R}$ is a constant, then we obtain the energy first introduced in \cite{Lei86}.  See also \cite{Tan96}.  We define the set, $K^{GL}$, of admissable maps associated with \eqref{GLhelfrich} by
\begin{align}
K^{GL}=\left\lbrace\begin{aligned}
(X,\phi)\in C^2(\Gamma_0;\mathbb{R}^3)\times H^1(\Gamma;\mathbb{R}):~&\mathcal{V}(X(\Gamma_0))=V_0
\\&\text{ and }\dashint_{\Gamma_0}\phi\:{\rm d}\mu=\alpha
\end{aligned}\right\rbrace
\end{align}
and consider the minimisation problem:
\begin{problem}\label{JulLipGL}
	Find $(X^*,\phi^*)\in K^{GL}$, such that
	$$\mathcal E_{JL}^{GL}(X^*,\phi^*)=\inf_{(X,\phi)\in K^{GL}}\mathcal E_{JL}^{GL}(X,\phi).$$
\end{problem} 
This minimisation problem for a diffuse interface energy can be understood, for small $\epsilon$, as an approximation to a minimisation problem of a sharp interface energy, \cite{EllSti10siam, EllSti10,EllSti12},  introduced in \cite{JulLip93,JulLip96}. A three dimensional phase field model was proposed  in \cite{WanDu08} in which the surface energy was also approximated by a bulk phase field energy.



 To write down the sharp interface energy we first let $\gamma_0$ denote a smooth curve on $\Gamma_0$ that carries line energy and divides $\Gamma_0$ into two non-empty regions $\Gamma_0^1$ and $\Gamma_0^2$, that is $\gamma_0=\partial \Gamma_0^1=\partial \Gamma_0^2$ and $\Gamma_0=\Gamma_0^1\cup\gamma_{0}\cup \Gamma_0^2$. 
Then given a mapping $X:\Gamma_0\to\Gamma$ we suppose that
$$\Gamma=\Gamma^1\cup\gamma\cup\Gamma^2$$
where $\Gamma^1$ and $\Gamma^2$ are defined by
\begin{align*}
\Gamma^i&=\left\{x\in \mathbb R^3~|~x=X(p),\,p\in \Gamma_0^i\right\}
\end{align*}
for $i\in\{1,2\}$ and $\gamma$ is defined by
\begin{align*}
\gamma&=\left\{x\in \mathbb R^3~|~x=X(p),\,p\in \gamma_0\right\}.
\end{align*}
The sharp interface energy is obtained from a generalised Canham-Helfrich energy functional by adding a line energy, \cite{JulLip93,JulLip96}, 
\begin{equation}
\mathcal E_{JL}(X):=\sum_{i=1}^2\int_{\Gamma^i_0}\left (\frac{\kappa_i}{2}(H-H_i^s)^2+\sigma\right){\rm d}\mu+\int_{\gamma} \sigma_\gamma \,{\rm d}\gamma.
\end{equation}
Here $\sigma_\gamma:=C_W b\in\mathbb{R}$ is a  line tension coefficient related to $b$ by the  constant $C_W$ depending only on the double well. We suppose that $\alpha\in (-1,1)$ is given. The corresponding sharp interface minimisation problem is then given by
\begin{problem}\label{JulLip}
Find $X^*\in K(X)$,
\begin{align}
K(X)=\left\lbrace\begin{aligned}
X\in C^1(\Gamma_0;\mathbb{R}^3):~&X|_{\Gamma_0^i}\in C^2(\Gamma_0^i;\mathbb{R}^3),\mathcal{V}(X(\Gamma_0))=V_0,\\
&\text{ and }|\Gamma^1|-|\Gamma^2|=\alpha|\Gamma|
\end{aligned}\right\rbrace
\end{align}
such that
$$\mathcal E_{JL}(X^*)=\inf_{X\in K}\mathcal E_{JL}(X).$$
\end{problem}
We anticipate that as $\epsilon \rightarrow 0$,  Problem \ref{JulLip} may  be obtained as a limit of Problem  \ref{JulLipGL}.

\section{Perturbation models of critical points}

A standard approach in biophysics is to consider  small deformations of flat surfaces written as graphs over the flat domain. This leads to variational problems for quadratic energy functionals involving the height of the graph together with appropriate boundary conditions and constraints modelling the interaction with individual particles, \cite{EllGraKor16}. 
These ideas were extended to critical points of the surface energy by \cite{EllFriHob17-a}. In particular the resulting quadratic energies were explicitly computed for the sphere and Clifford torus. In the following,  after considering flat domains and the Monge gauge, we describe the calculations of \cite{EllFriHob17-a}.

 \subsection{Flat domains and the Monge gauge}\label{subSec:SmallDefPlane}

In the Monge gauge, one first assumes that the surface is nearly flat in the sense that   it can be written as a graph
\begin{align}
    \label{eq:monge}
    {\Gamma}=\{(x_{1},x_{2},u(x_{1},x_{2})) | (x_1,x_2)\in \Omega\}
\end{align} 
over a two-dimensional reference domain $\Omega\subset \mathbb R^2$.
Then, the mean curvature $H$ and the Gaussian
curvature $K$ of the membrane $\Gamma$ are given by
\begin{equation*}
    H = -\nabla \cdot \frac{\nabla u}{(1+|\nabla u|^{2})^{1/2}}, ~~
    K = \left(\frac{\partial^{2}u}{\partial x_1^{2}}\frac{\partial^{2}u}{\partial x_2^{2}} 
    - \left(\frac{\partial^{2}u}{\partial x_1\partial x_2}\right)^2\right) \left /  {(1+|\nabla u|^{2})^{1/2}}\right. .
\end{equation*} 

A common approach to derive an approximate model is to assume 
that the displacement of the membrane from the $(x_1,x_2)$ plane produced by the particles is small, i.e. $|\nabla u|\ll 1$. 
This leads to the well known  geometric linearisation
\begin{equation} \label{eq:GEOLIN}
 (1+|\nabla u|^{2})^{1/2} \; \rightsquigarrow \; 1 +{\textstyle \frac{1}{2}}|\nabla u|^{2},\quad
    H \;\rightsquigarrow \; -\Delta u,\quad
    K\; \rightsquigarrow \; \frac{\partial^{2}u}{\partial x_1^{2}}\frac{\partial^{2}u}{\partial x_2^{2}} - \left(\frac{\partial^{2}u}{\partial x_1\partial x_2}\right)^2,
\end{equation}
modelling perturbations from a flat surface
yielding, up to a constant term, the quadratic energy 
\begin{equation*}\label{energy_monge}
    \mathcal E(u) = \int_{\Omega} \left (\textstyle \frac{1 }{2}\kappa(\Delta u)^2 
    + \kappa_G \left(\frac{\partial^{2}u}{\partial x_1^{2}}\frac{\partial^{2}u}{\partial x_2^{2}} - 
    \left(\frac{\partial^{2}u}{\partial x_1\partial x_2}\right)^2\right)  + \frac{1 }{2} \sigma |\nabla u|^2\right ) \; {\rm d}x.
\end{equation*}
Finally, ignoring Gaussian curvature as discussed above, we obtain 
a quadratic approximation of the energy taking the form
\begin{equation}\label{eq:monge_energy}
    \mathcal E(u)= 
    \int_\Omega \left ({\textstyle \frac{1}{2}} \kappa(\Delta u)^2  + {\textstyle \frac{1}{2}} \sigma|\nabla u|^{2} \right )\; {\rm d}x .
\end{equation}
The deformation of the membrane from the flat configuration may be  in response  to small forces or deformation constraints. The minimization of the resulting energy leads  to fourth order  equations involving the biharmonic operator. These are considerably simpler than the fourth  order nonlinear equations  arising from minimising the original energy.

 \begin{remark}
The energy \eqref{eq:monge_energy} with appropriate perturbations  is systematically studied in \cite{EllGraKor16} in the case of point constraints arising from attachments  of the membrane to proteins modelled as points. A phase field perturbation is studied in \cite{FonHayLeo16}. \end{remark}

\subsection{First and second variations of geometric quantities on spheres}\label{subSec:FirstSecVar}
We now outline the derivation of the approximating quadratic energy functional in the case 
that the membrane is a small displacement of a sphere. Following \cite{EllFriHob17-a}, the idea is to take advantage of the fact that  the sphere is a critical point of the Willmore energy. 

\subsubsection{Variations of functionals}

It is convenient to  define the following Lagrangian
\begin{equation}
\mathcal{L}(\Gamma,\lambda):= \kappa \mathcal W(\Gamma) +\sigma \mathcal A(\Gamma)  +\lambda\left(\mathcal V(\Gamma)-V_0\right),
\end{equation}
where $\lambda$ is the Lagrange multiplier associated with the volume constraint.

\begin{remark}
$\lambda$  can be interpreted as a hydrostatic pressure.  
\end{remark}

We consider graph like  surfaces of the form,
\begin{equation}\label{eq:perturbedSurface}
\Gamma_\rho(u) := \{ x + \rho u(x)\nu(x) : x \in \Gamma_0\}
\end{equation}
where $u \in C^2(\Gamma_0)$, $\nu$ is the outward pointing unit normal to $\Gamma_0=\mathbb{S}^2(R)$ and $\rho$ is a constant. 
We will assume $\rho\ll 1$ and hence only consider small perturbations.
To compute the required derivatives of a functional  $J$ we use the following formulae relating them to variations of $J$.
\begin{align*}
 J'(\Gamma_0)[u\nu] &:= \frac{\mathrm{d}J(\Gamma_\rho)}{\mathrm{d}\rho}\bigg|_{\rho =0}  \\
  J''(\Gamma_0)[u\nu,u\nu] & := \frac{\mathrm{d^2}J(\Gamma_\rho)}{\mathrm{d}\rho^2}\bigg|_{\rho =0}.
\end{align*}
For our case, the second variation will coincide with the variation of the first variation (see Remark 3.2 in \cite{EllFriHob17-a}) given by 
\begin{eqnarray}\label{varOfVar}
J^{\prime\prime}(\Gamma_0)[g\nu,u\nu]:=\left.\frac{\mathrm{d}J^\prime(\Gamma_\rho)[g\nu_\rho]}{\mathrm{d}\rho}\right|_{\rho=0},
\end{eqnarray}
where $\nu_\rho$ is the outward unit normal to $\Gamma_\rho$. Hence we will use \eqref{varOfVar} below to calculate the second variation.

The first and second variations of the  Willmore functional on a closed surface $\Gamma_0$ are
\begin{equation}
\label{Willmore_functional_1var}
	\mathcal W'(\Gamma_0)[u\nu] = \int_{\Gamma_0} \left ( -\Delta_{\Gamma_0} H  -  H|\mathcal{H}|^2 +\frac{1}{2}H^3 \right )u \;{\rm d}\Gamma_0.
\end{equation}

\begin{align}\label{willmore2varIBP}
\begin{split}
& \mathcal W''({\Gamma_0})[u\nu, g\nu]  = \int_{\Gamma_0}  (\Delta_{\Gamma_0} u + |\mathcal{H}|^2u)(\Delta_{\Gamma_0} g + |\mathcal{H}|^2g) \;{\rm d}{\Gamma_0} \\ 
& \quad  +\int_{\Gamma_0} 2H \mathcal{H}:(g\nabla_{\Gamma_0} \nabla_{\Gamma_0} u + u\nabla_{\Gamma_0} \nabla_{\Gamma_0} g )\;{\rm d}{\Gamma_0} \\
& \quad +\int_{\Gamma_0} \left (2H \mathcal{H}\nabla_{\Gamma_0} u \cdot \nabla_{\Gamma_0} g - \frac{3}{2} H^2 \nabla_{\Gamma_0} u \cdot \nabla_{\Gamma_0} g 
-\frac{3}{2}H^2 ( u \Delta_{\Gamma_0} g + g \Delta_{\Gamma_0} u )\right ){\rm d}{\Gamma_0}\\
& \quad +\int_{\Gamma_0}  ( 2HTr(\mathcal{H}^3) -\frac{5}{2}H^2|\mathcal{H}|^2 + \frac{1}{2}H^4) gu  \;{\rm d}{\Gamma_0} 
\end{split}
\end{align}

The first and second variation of the area functional $\mathcal A({\Gamma_0})$ are given by 
\begin{align}
	&	\mathcal A'({\Gamma_0})[u\nu] = \int_{\Gamma_0} u H \; {\rm d}{\Gamma_0}, \\
	&	\mathcal A''({\Gamma_0})[u\nu, g\nu] = \int_{\Gamma_0}  \nabla_{\Gamma_0} u \cdot \nabla_{\Gamma_0} g + (H^2 - |\mathcal{H}|^2) u g  \; {\rm d}{\Gamma_0}.
\end{align}
The first and second variation of the volume functional $V({\Gamma_0})$ are given by
\begin{align}
& \mathcal V'({\Gamma_0})[u\nu] = \int_{\Gamma_0} u \;{\rm d}{\Gamma_0}, \label{volume1var}\\
& \mathcal V''({\Gamma_0})[u\nu,g\nu] = \int_{\Gamma_0} Hgu \;{\rm d}{\Gamma_0}.
\end{align}
 Writing $\lambda_\rho=\lambda_0+\rho \mu $ we may expand
\begin{equation} \label{L_taylor_series}\begin{split}
\mathcal{L}(\Gamma_\rho,\lambda_\rho) & = \mathcal{L}(\Gamma_0,\lambda_0) + \rho\frac{\mathrm{d}\mathcal{L}(\Gamma_\rho,\lambda_\rho)}{\mathrm{d}\rho}\bigg|_{\rho=0}   \\
& + \frac{\rho^2}{2}  \frac{\mathrm{d}^2\mathcal{L} (\Gamma_\rho,\lambda_\rho)}{\mathrm{d}\rho^2}\bigg|_{\rho=0}  
+ O(\rho^3).
\end{split}
\end{equation}

We follow  \cite{EllFriHob17-a}  where the idea is to exploit the vanishing of the first variation of the Lagrangian at $(\Gamma_0,\lambda_0)$ and  neglect the  $O(\rho^3)$ term yielding an approximation based on a quadratic Lagrangian.

\begin{defn}
For the surface $\Gamma_0 \subset \mathbb{R}^{3}$ with associated Lagrange multiplier $\lambda_0 \in \mathbb{R}$, the quadratic Lagrangian $L:H^2(\Gamma_0) \times \mathbb{R} \rightarrow \mathbb{R}$ is given by
\begin{align*}
L(u,\mu):&=\frac{1}{2} \frac{\mathrm{d}^2\mathcal{L}(\Gamma_\rho(u),\lambda(\rho))}{\mathrm{d}\rho^2}\bigg|_{\rho=0} \\
& = \frac{1}{2}\mathcal{E}''(\Gamma_0)[u\nu,u\nu] + \frac{1}{2}\lambda_0 \mathcal{V}''(\Gamma_0)[u\nu,u\nu] + \mu \mathcal{V}'(\Gamma_0)[u\nu]. 
\end{align*}
\end{defn}

\subsubsection{Variations of functionals on sphere}

 Since a sphere has zero first variation of the Willmore functional,  the first variation of the Lagrangian with respect  to $\Gamma$ 
 evaluated on $\mathbb{S}^2(R)$ in the normal direction is
$$\sigma\frac{2}{R}+\lambda.$$
 Setting $\lambda_0=-2\sigma/R$ we observe  that the first variation vanishes when $\lambda=\lambda_0$ so
  $(\Gamma_0,\lambda_0)$ is a critical point of the Lagrangian  $\mathcal{L}(\cdot,\cdot)$ with $V_0=V(R)$.

The second variations on $\Gamma_0=\mathbb{S}^2(R)$ may be calculated to be

\begin{align*}
& \mathcal W''(\Gamma_0)[u\nu,g\nu]=\int_{\Gamma_0}\ \left (\Delta_{\Gamma_0} u\Delta_{\Gamma_0} g-\frac{2}{R^2}\nabla_{\Gamma_0}u\cdot\nabla_{\Gamma_0} g\right ) {\rm d}\Gamma_0, \\
&  \mathcal A''(\Gamma_0)[u\nu,g\nu]=\int_{\Gamma_0}\left (\nabla_{\Gamma_0} u\cdot\nabla_{\Gamma_0} g+\frac{2}{R^2}ug \right ){\rm d}\Gamma_0, \\ 
& \mathcal V''(\Gamma_0)[u\nu_0,g\nu_0] =\frac{2}{R}\int_{\Gamma_0}ug {\rm d}\Gamma_0. 
\end{align*}
It follows that

$$L(u,\mu)=\frac{\kappa}{2}\mathcal{W}''(\Gamma_0)[u\nu,u\nu]+\frac{\sigma}{2}\mathcal{A}''(\Gamma_0)[u\nu,u\nu]-\frac{\sigma}{R}\mathcal{V}''(\Gamma_0)[u\nu,u\nu]+\mu\mathcal V'(\Gamma_0)[u\nu]$$
yielding

\begin{equation}
L(u,\mu)=\int_{\Gamma_0}\left ( \frac{\kappa}{2}|\Delta_{\Gamma_0} u|^2+\frac{1}{2}\left(\sigma-\frac{2\kappa}{R^2}\right)|\nabla_{\Gamma_0}u|^2-\frac{\sigma}{R^2}u^2+\mu u\right ){\rm d}\Gamma_0.
\end{equation}
so that
\begin{equation}
L(u,\mu)=\frac{1}{2}a(u,u)+\mu (u,1)_{L^2(\Gamma_0)}.
\end{equation}

Here  we have defined  the bilinear form associated to the quadratic part of $L(\cdot,\cdot)$, $a(\cdot,\cdot)\colon H^2(\Gamma_0) \times H^2(\Gamma_0) \to \mathbb{R}$ as follows,
\begin{equation}\label{eq:sphereBilinear}
a(u,v) = \int_{\Gamma_0} \left ( \kappa \Delta_{\Gamma_0} u \Delta_{\Gamma_0} v + \left(\sigma - \frac{2\kappa}{R^2}\right)\nabla_{\Gamma_0} u \cdot \nabla_{\Gamma_0}v - \frac{2\sigma}{R^2}uv\right ){\rm d}\Gamma_0.
\end{equation}
It can be shown that $a(\cdot,\cdot)$ is coercive over the space $\{1,\nu_1,\nu_2,\nu_3\}^\perp$ \cite{EllFriHob17-a}, where $\perp$ is meant in the sense of $H^2(\Gamma_0)$. However, since these functions are positive eigenfunctions of $-\Delta_{\Gamma_0}$ it is equivalent to take the $L^2(\Gamma_0)$ orthogonality. Note on  $\mathbb{S}^2(R)$ that  $\nu=x/R$. 
The orthogonality constraint
$$(1,u)=0$$  may be interpreted as enforcing the volume constraint and 
$$(x_i,u)=0,~i=1,2,3$$
as ensuring the surface has a fixed center of mass.
\begin{remark}
The Monge-Gauge approximation  \eqref{eq:monge_energy} for near flat membranes is recovered by formally taking $R\to \infty$ in  \eqref{eq:sphereBilinear}.
 \end{remark}
\subsubsection{Perturbed energy functional}\label{PerturbAbs}

We wish to consider critical points  of the perturbed Lagrangian
\begin{equation}\label{eq:NonLin}
\mathcal{L}_\rho (\Gamma_\rho(u),\lambda_\rho) = \mathcal{L}(\Gamma_\rho(u),\lambda_\rho)+\mathcal{F}_\rho(\Gamma_\rho(u)),
\end{equation}
 Here $\mathcal{F}_\rho(\cdot)$ is a perturbative  energy whose precise form depends  on the type of deformation being considered. Here we supose that it may be expanded as
  \begin{equation}
  \mathcal F_\rho(\Gamma_\rho)=\rho\mathcal F_0(\Gamma_0)+\rho^2 \mathcal F_1(\Gamma_0)[u\nu]+\mathcal{O}(\rho^3).
 \end{equation} 
 We write the Taylor expansion of \eqref{eq:NonLin} with $\Gamma_\rho(u)$ set as in (\ref{eq:perturbedSurface}) and $\lambda_\rho=\lambda_0+\mu\rho$  up to third order in $\rho$ as
\begin{align*}
\mathcal{L}_\rho (\Gamma_\rho(u),\lambda_\rho) & = \mathcal{L}(\Gamma_0,\lambda_0)\\
 & +\rho \left ( \mathcal E'(\Gamma_0)[u\nu] +\lambda_0\mathcal V'(\Gamma_0)[u\nu]+\mu(\mathcal V(\Gamma_0)-V_0) +\mathcal{F}_0(\Gamma_0) \right ) \\
& +\frac{ \rho^2}{2} \left (\mathcal E''(\Gamma_0)[u\nu,u\nu]+\frac{\lambda_0}{2}\mathcal V''(\Gamma_0)[u\nu,u\nu]+2\mu\mathcal V'[u\nu] +2\mathcal F_1(\Gamma_0)[u\nu]\right )\\
& + \mathcal{O}(\rho^3).
\end{align*}
 Evauating on the critical point  $(\mathbb{S}^2(R),\lambda_0=\frac{-2\sigma}{R})$  of  $\mathcal L(\cdot,\cdot)$ we find
 \begin{align*}
 \mathcal{L}_\rho (\Gamma_\rho(u),\lambda(\rho)) & = \mathcal{L}(\Gamma_0,\lambda_0)+\rho\mathcal F_0(\Gamma_0)\\
 & +\frac{ \rho^2}{2}\left(L(u,\mu)+2\mathcal F_1(\Gamma_0)[u\nu]\right )+ \mathcal{O}(\rho^3).
 \end{align*}
 
 It follows that up to third order in $\rho$ a critical point of  \eqref{eq:NonLin} may be approximated by finding a critical point $(u,\mu)$ of
 \begin{equation}\label{quadLag}
 \mathcal J(u,\mu):=L(u,\mu)+2 F(u)
 \end{equation}
 where $F(u):=\mathcal F_1(\Gamma_0)[u\nu]$.


\section{Point constraints}\label{Points}

Here we further explore problems where we wish to minimise the energy \eqref{Willmore+tension} subject to point constraints as described in Problems   \ref{soft} and  \ref{hard} . From the physical perspective this models the interaction of proteins with a biomembrane.
To include such mechanisms, we wish to allow for groups of point constraints to move rigidly.
A group of point constraints may be considered to be representative of a large biological structure, for example a BAR protein \cite{HenKenFor07}, which  are believed to attach to the membrane at a finite number of points.
We consider $N$ such groupings of identical particles, each of which has $L$ points of attachment.  
Following \cite{Kie18}, for $q \in \mathbb{R}^6$, we write
\begin{equation}\label{eq:6DoF}
\phi(q, x) := R_x(q_1)R_y(q_2)R_z(q_3)x + (q_4,q_5,q_6)^T
\end{equation}
to be the map describing the rigid transformation, where $R_x$, $R_y$ and $R_z$ are the rotations in the $x$, $y$ and $z$ axes respectively.
This then allows us to pose the Problems   \ref{soft} and  \ref{hard} in the setting of the surface being a graph with rigid transformations of particles.

\begin{problem}
	Find $u$, $\{q^i\}_{i=1}^N$  minimising
	\begin{equation}\label{eq:PenaltyEquation}
		\int_{\Gamma_\rho(u)} \left( \frac{\kappa}{2} H^2 + \sigma \right) {\rm d} \Gamma_\rho(u) + \frac{1}{2\delta}\sum_{i=1}^N\sum_{j=1}^L d(\Gamma_\rho(u),\phi(q^i,z_j) )^2
	\end{equation}
	subject to $\mathcal{V}(\Gamma_{\rho}(u)) = V_0$.
\end{problem}
\begin{problem}
	Find $u$, $\{q^i\}_{i=1}^N$ minimising \eqref{Willmore+tension} subject to $\phi(q^i,z_j) \in \Gamma_\rho(u)$ for $j=1,...,L$, $i=1,...,N$ and $\mathcal{V}(\Gamma_{\rho}(u)) = V_0$.
\end{problem}

\begin{remark}
Something we have neglected to consider is that the protein is not only a collection of points which are attached to the membrane,
but also may be  inside or outside the interior of the surface. This may be modelled with  inequality constraints.

\end{remark}
We see that these problems are highly non-linear and should they be solvable, numerical methods will be challenging.
We wish to exploit $\rho$ as a small parameter writing 
$$\mathcal F_\rho(\Gamma_\rho)=   \frac{1}{2\delta}\sum_{i=1}^N\sum_{j=1}^L d(\Gamma_\rho(u),\phi(q^i,z_j(\rho)) )^2$$
in order to simplify the problem as discussed in sub-section \ref{PerturbAbs}. Here as  we assume the deformations are small, we  write the points, $z_j$, as a graph over $\Gamma_0$, specifically,
\begin{equation}
	z_j = p_j + \rho Z_j \nu(p_j)
\end{equation}
for each $j=1,...,L$, where $Z_j\in \mathbb{R}$ and the $p_j\in \Gamma_0$ are distinct. 

We now limit ourselves to a single group of constraints ($N=1$) and we fix $q$ so that we  discuss only the geometric minimisation for a fixed location of a single particle.

As in \cite{EllFriHob17-a}, we calculate
\begin{align*}
d(\Gamma_\rho( u),p_jh +\rho Z_j \nu (p_j) )|_{\rho = 0}
=
d(\Gamma_0, p_j) = 0,
\end{align*}
with first derivative term,
\begin{align*}
\frac{{\rm d}}{{\rm d} \rho} d(\Gamma_\rho (u),p_j +\rho Z_j \nu (p_j) )|_{\rho = 0}
=&
\nabla d (\Gamma_0,p_j) \cdot \frac{{\rm d}}{{\rm d}\rho} (p_j + \rho Z_j \nu(p_j) )|_{\rho=0}
\\ &+ \dot{\partial} d(\Gamma_\rho (u),p_j)|_{\rho=0}
\\
=& Z_j - u(p_j),
\end{align*}
which results in
\begin{equation}
	d(\Gamma_\rho (u),p_j +\rho Z_j \nu (p_j) )^2
	=
	\rho^2 (Z_j -u(p_j))^2 + \mathcal{O}(\rho^3)
\end{equation}
and
\begin{equation}
F(u)=\frac{1}{2\delta}\sum_{j=1}^L |u(p_j)-Z_j|^2.
\end{equation}


The Lagrangian (\ref{quadLag}) becomes
\begin{equation}
\mathcal J(u,\mu)=\frac{1}{2}a(u,u)+\mu(1,u)+\frac{1}{2\delta}\sum_{j=1}^L|u(p_j)-Z_j|^2.
\end{equation}
 By considering
 \[ U_0:=\{ v \in H^2({\Gamma_0})~|~(1,v)=0\}\]
we may now consider finding a minimiser of the quadratic functional \begin{equation}
\mathcal Q(u):=\frac{1}{2}a(u,u)+\frac{1}{2\delta}\sum_{j=1}^L|u(p_j)-Z_j|^2
\end{equation} over $U_0$. In \cite{Hob16} it is seen that the largest subspace of  $U_0$ where  $a$ is coercive is
\[
	U_{\nu}:=\{ v \in H^2({\Gamma_0}) : v\in (\text{Sp}\{1,\nu_1,\nu_2,\nu_3\})^\perp \}.
\]
In the case of $L\geq 4$ and there is $j_1,j_2,j_3,j_4$ with $\{p_{j_k}\}_{k=1}^4$ non-coplanar, it may be seen that $Q$ is coercive over $U_0$. 
We now have the following existence and uniqueness results, which follow from the Lax-Milgram theorem.
\begin{theorem}\label{thm:penalty}
Given distinct  $\{p_i\}_{i=1}^L\subset \Gamma_0$ and $\{Z_i\}_{i=1}^L\subset \mathbb{R}$, with penalty parameter $\delta>0$, there is a unique $u^\delta\in U_\nu$ such that for any $v \in U_\nu$,
\[
	a(u^\delta,v) + \frac{1}{\delta} \sum_{j=1}^L u^\delta(p_j) v(p_j) = \frac{1}{\delta} \sum_{j=1}^L Z_j v(p_j).
\]
\end{theorem}

\begin{theorem}\label{thm:hardConstraint}
Given disjoint  $\{p_i\}_{i=1}^L\subset \Gamma_0$ , and $\{Z_i\}_{i=1}^L\subset \mathbb{R}$, there is a unique $u\in U_\nu$ such that $u(p_i) = Z_i$ for $i=1,...,L$ and for any $v \in U_\nu$ with $v(p_i) = 0$ for $i=1,...,L$, it holds
\[
	a(u,v) = 0.
\]
\end{theorem}

Viewing $u^\delta$ as a penalised approximation to $u$, an interesting question to now ask is, as $\delta \to 0^+$, does $u^\delta \to u$?
If so, in what sense does it converge (i.e. what topology)?
If it converges in a norm, can we show a reasonable error estimate?
These questions lead to the following convergence result.
\begin{theorem}
	Let $u^\delta \in U_\nu$ be as in Theorem \ref{thm:penalty}, let $u\in U_\nu$ be as in Theorem \ref{thm:hardConstraint}, with the constraints $\{p_j\}_{j=1}^L\subset \Gamma_0$ and $\{Z_j\}_{j=1}^L\subset \mathbb{R}$.
	Then there is $C>0$ independent of $\delta$ such that
	\[
		\|u-u^\delta\|_{H^2(\Gamma_0)} \leq C \sqrt{\delta}.
	\]
\end{theorem}
We have previously fixed $q$ for convenience since finding a minimiser for fixed $q$ is a standard variational problem, whereas the problem for minimising $q$ is not variational.
The problem for finding a minimising $q$ is developed as a gradient descent in \cite{Kie18}, which deals with a near flat membrane and also that of a near cylindrical membrane.

\begin{figure}[h]
\includegraphics[width=.5\linewidth]{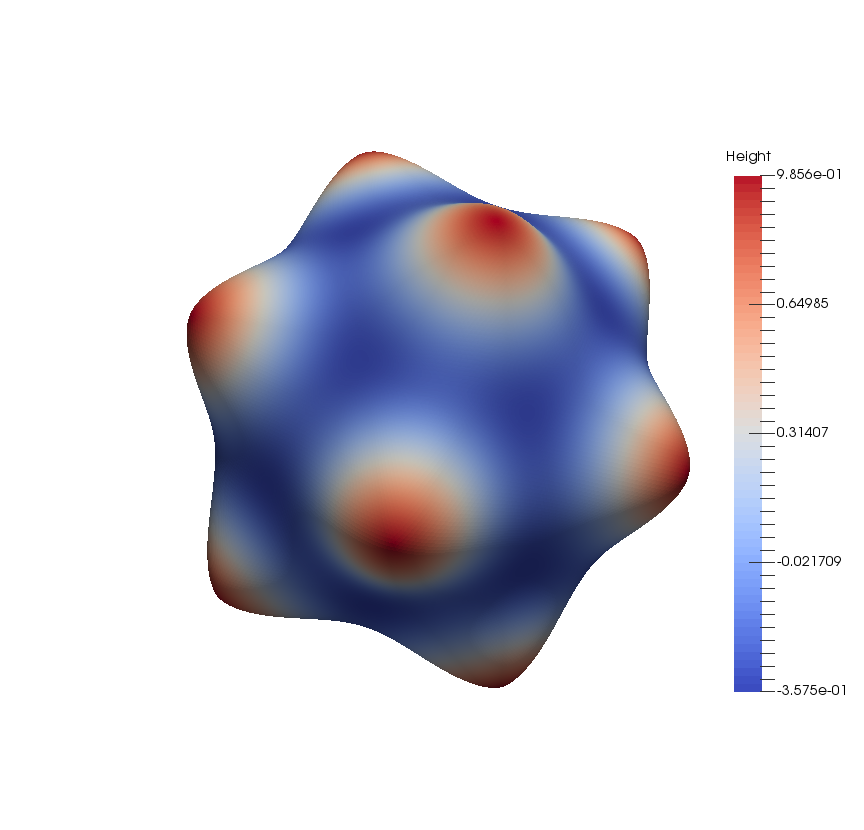}
\caption{Example numerical simulation with $\delta = 10^{-4}$ and 12 points  $X_i$ positioned on the vertices of a regular icosahedron with heights $1$, $\rho = 0.2$.}
\label{fig:icos}\end{figure}

\begin{figure}[h]
\includegraphics[width=.45\linewidth]{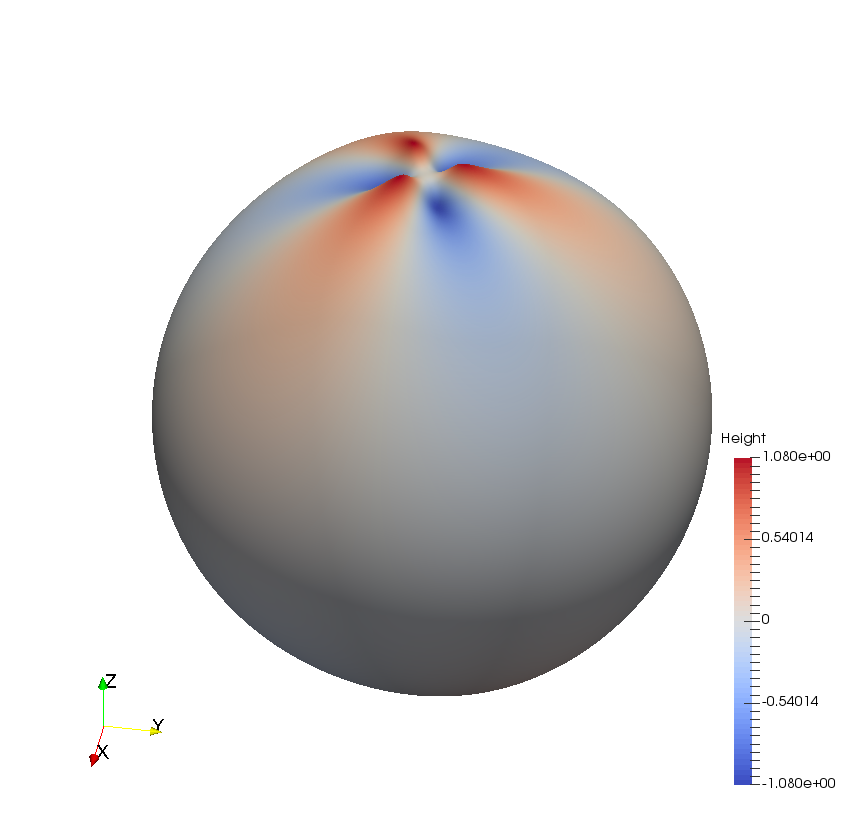}
\includegraphics[width=.45\linewidth]{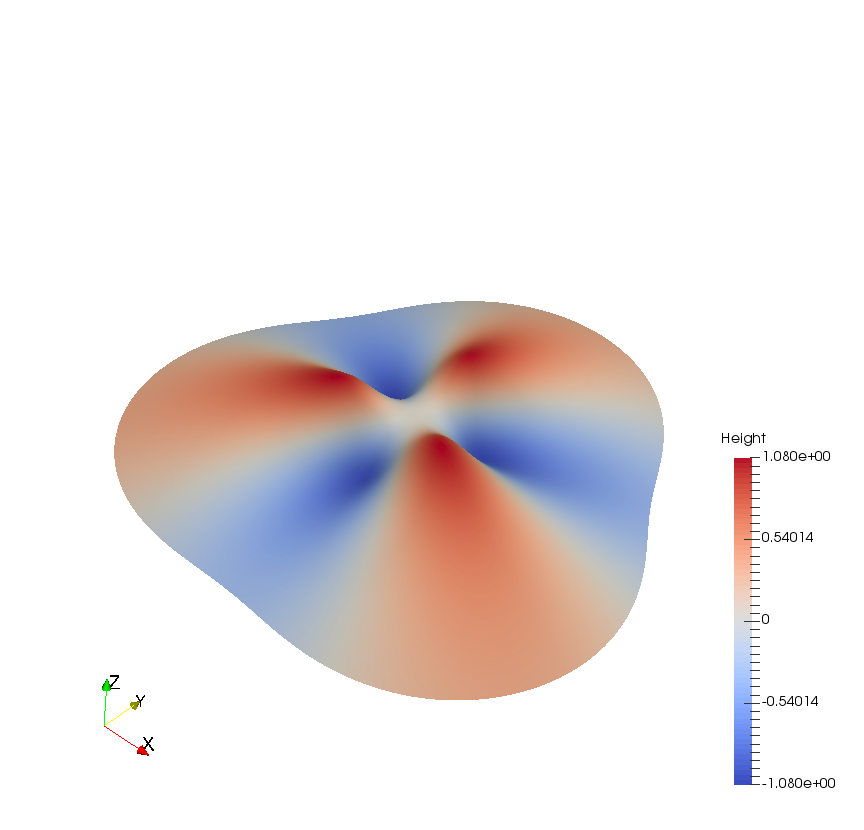}
\caption{Example simulation with $\delta = 10^{-6}$,  two particles   positioned around each of the North and South poles with 6 points of attachment, placed  with heights alternating between $-1$ and $1$, antisymmetrically, $\rho = 0.05$.
Left: whole sphere; Right: zoom in of North pole}
\label{fig:NandS}
\end{figure}

\begin{figure}[h]
\includegraphics[width=0.8\linewidth]{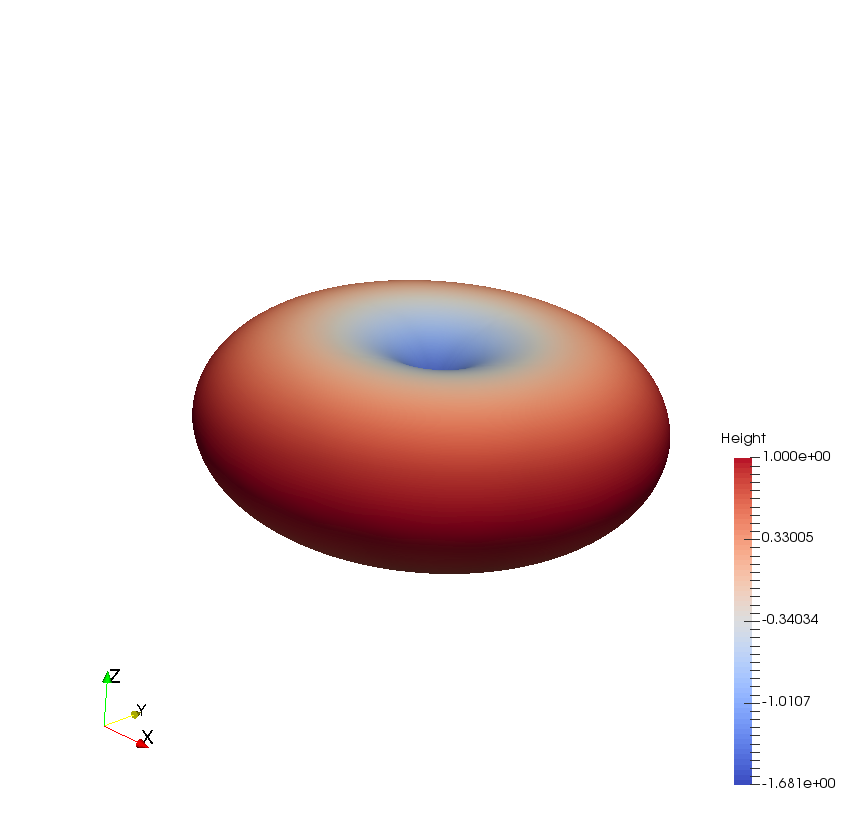}
\caption{Example simulation with $\delta = 10^{-6}$, with 10 points of attachment  equally spaced on the equator, with height 1, $\rho = 0.6$.}
\label{fig:redblood}\end{figure}

The methods for numerically solving the penalty problems are based on splitting the fourth order problem into coupled second order problems \cite{EllFriHob19}.
The methods for extending this to the case of the hard constraint problem are to be found in a work in preparation by Elliott and Herbert.
The simulations illustrated in Figures \ref{fig:icos}, \ref{fig:NandS} and \ref{fig:redblood}  are implemented using the DUNE framework \cite{
BasBlaDed08-b},
 in particular the ALU-Grid module \cite{AlkDedKlo16}.

\section{Phase field interaction}\label{PhaseField}
The second perturbation model we review takes into account phase separation between two differing lipid compositions within the membrane.   Phase separation has been observed in artificial membranes such as in giant unilamellar vesicles (GUVs) or supported lipid bilayers (SLBs). Experimental results involving GUVs and SLBs suggest the curvature of the biological membrane plays an important role in the formation of these domains \cite{BauDasWeb05,ParYuChe06,RinFonGio18,RozKaiYos05,BauHesWeb03}.  Lipid rafts  are small (10-200nm) domains which compartmentalise cellular functions such as signal transduction, protein sorting and host-pathogen interactions \cite{Pik06,SezLevMay17}.

Here we introduce a concentration $\phi$ taking the values $\pm 1$ in each of two phases and consider 
 the energy
\begin{align}
\label{Phase-Energy}
\mathcal E(\Gamma,\phi)=&\kappa\mathcal W(\Gamma)+\sigma\mathcal A(\Gamma) - \kappa\Lambda \mathcal I(\Gamma,\phi) 
+ b \mathcal P(\Gamma,\phi)
\end{align} 

subject to the following  concentration  and volume constraints
\begin{align}
&\mathcal R(\Gamma,\phi):=\dashint_{\Gamma}\phi=\alpha,&&\mathcal{V}(\Gamma)=V(R).
\end{align}
Here $\mathcal{I}(\cdot,\cdot)$ and $\mathcal{P}(\cdot,\cdot)$ are given by
\begin{align}
\mathcal I(\Gamma,\phi):= &\int_\Gamma H\phi {\rm d}\Gamma,\\
\mathcal P (\Gamma,\phi): =&
\int_\Gamma \left (\frac{\epsilon}{2}|\nabla_\Gamma\phi|^2+\frac{1}{\epsilon}f(\phi)\right ) {\rm d}\Gamma
\end{align}
and where
\begin{align}
f(\phi):=&W(\phi)+\frac{\epsilon\kappa\Lambda^2\phi^2}{2b}.
\end{align}

This leads to the consideration of the Lagrangian
\begin{equation}
\mathcal L (\Gamma,\phi,\lambda,\eta):=\mathcal E(\Gamma,\phi)+\lambda (\mathcal V(\Gamma)-V_0)+\eta(\mathcal R(\Gamma,\phi)-\alpha)
\end{equation}

 where $\lambda$ and $\eta$ are the Lagrange multipliers for the constraints.

The  energy (\ref{Phase-Energy})   may be rewritten to be in the form  \eqref{GLhelfrich} in which  the protein molecules induce a spontaneous curvature of the form
\begin{displaymath}
H^s(\phi)=\Lambda\phi.
\end{displaymath}
This model was  proposed in \cite{Lei86}. Similar models have been considered in \cite{EllSti10,EllSti10siam,  EllSti12, HeaDha17,WanDu08} where phase-dependent bending rigidities were considered.




We now suppose that the spontaneous curvature is small leading to a perturbation of the underlying equilibrium.  It is  reasonable to expect the deformations to be of the same order, and hence we scale both of these terms by $\rho$.  In \cite{KuzAkiChi05} it was calculated that the line tension for lipid rafts should depend quadratically on the spontaneous curvature, so we  scale this term by $\rho^2$. 
This leads to  a perturbation model for which we have the scaling
\begin{align*}
&\Gamma\rightsquigarrow \Gamma_\rho,&&\Lambda\rightsquigarrow\rho\Lambda,&& b\rightsquigarrow\rho^2 b.
\end{align*}

Carrying out calculations similar to those given earlier leads to an energy
\begin{equation}\label{Raft-Energy}
E(u,\phi):=\frac{1}{2}a(u,u)+I(u,\phi)+G(\phi)
\end{equation}
where
\begin{align}
I(u,\phi):=&\int_{\Gamma_0} \left(\kappa\Lambda\phi\Delta_{\Gamma_0}u+\frac{2\kappa\Lambda u\phi}{R^2}\right ){\rm d}\Gamma_0,\\
G(\phi):=&b\int_{\Gamma_0} \left(\frac{\epsilon}{2}|\nabla_{\Gamma_0}\phi|^2+
\frac{1}{\epsilon}f(\phi)\right ){\rm d}\Gamma_0.
\end{align}
The perturbed equations for the volume and mean concentration constraints become
\begin{equation}\label{constraint-phi/u}
\int_{\Gamma_0}u=0\qquad\text{ and }\qquad\dashint_{\Gamma_0}\phi=\alpha.
\end{equation}

Formally considering $R\to\infty$ we  obtain the approximation in the case $\Gamma_0$ is a plane. This is in agreement with the linearisation outlined in \cite{EllGraKor16,FonHayLeo16} and the energy used in \cite{AndKawKaw92,KawAndKawTan1993, Lei86,LeiAnd87,SeuAnd95}.

 Observe that for components of the normal $\nu_i$ we have that
\begin{equation}
E(u+\tau\nu_i,\phi)=E(u,\phi) \text{ for all } \tau\in\mathbb{R} \text{ and } i\in\{1,2,3\}.
\end{equation} 
This is reasonable, since this corresponds to order $\rho$ translations of the centre of mass, for which we would expect the energy to be invariant. \\
\\
We can therefore further constrain (\ref{Raft-Energy}) by requiring that
\begin{equation}\label{normal-constraint}
\int_{\Gamma_0}u\nu_i=0\qquad\text{ for }i\in\{1,2,3\}
\end{equation}
to remove this redundancy. This has the added benefit that $E(\cdot,\cdot)$ is coercive over this subspace. By considering $h=\tau\nu_1$ and $\phi=\alpha$, it follows that \eqref{normal-constraint} is neccessary for coercivity since 
\begin{displaymath}
\lim_{\tau\to\infty}\left(\|\tau\nu_1\|_{H^2(\Gamma_0)}+\|\alpha\|_{H^1(\Gamma_0)}\right)=\infty
\end{displaymath}
but
\begin{displaymath}
{E}(\tau\nu_1,\alpha)=\frac{b}{\epsilon}\int_{\Gamma_0}f(\nu_1)\leq C
\end{displaymath}
where $C$ is some constant independent of $\tau$. In fact by defining
\begin{displaymath}
\mathcal{K}:=\left\{(u,\phi)\in H^2(\Gamma_0)\times H^1(\Gamma_0):u\in \text{Sp}\left\{1,\nu_1,\nu_2,\nu_3\right\}^\perp \text{ and }\dashint_{\Gamma_0}\phi=\alpha\right\}.
\end{displaymath}
it follows that this is sufficient since we can prove there exists positive constants $C_1$ and $C_2$ such that the energy functional ${E}(\cdot,\cdot):\mathcal{K}\to\mathbb{R}$ satisifies the following inequality
	\begin{displaymath}
	{E}(u,\phi)\geq C_1\left(\|\phi\|^2_{H^1(\Gamma_0)}+\|u\|^2_{H^2(\Gamma_0)}\right)-C_2\qquad\forall (\phi,u)\in\mathcal{K}.
	\end{displaymath}
Hence it can be shown using the direct method of calculus of variations  that there exists a solution to the following problem.
\begin{problem}
	\label{subspaceProblem}
	Find $({u^*},{\phi^*})\in\mathcal{K}$ such that
	\begin{displaymath}
	{E}({u^*},{\phi^*})=\inf_{(u,\phi)\in\mathcal{K}}{E}(u,\phi).
	\end{displaymath}
\end{problem}

Finally, in order to calculate critical points of $E(\cdot,\cdot)$ over $\mathcal{K}$ it is convenient to investigate these numerically using a conserved $L^2$-gradient flow. This yields the  following system  
\begin{align}
\label{AC PDE}
\begin{cases}
\alpha_1\phi_t+\frac{b}{\epsilon}f^\prime(\phi)-b\epsilon\Delta_{\Gamma_0}\phi+\kappa\Lambda\Delta_{\Gamma_0} u+\frac{2\kappa\Lambda u}{R^2}+\lambda_\phi=0\\
\alpha_2 u_t-\left(\sigma-\frac{2\kappa}{R^2}\right)\Delta_{\Gamma_0} u+\kappa\Delta_{\Gamma_0}^2 u-\frac{2\sigma u}{R^2}+\kappa\Lambda\Delta_{\Gamma_0}\phi+\frac{2\kappa\Lambda \phi}{R^2}+\lambda_u=0\\
\end{cases}
\end{align}
subject to some initial conditions which satisfy the constraints \eqref{constraint-phi/u} and \eqref{normal-constraint}. These constraints are preserved by the evolution when values of the Lagrange multipliers corresponding to the constraints \eqref{constraint-phi/u} have the following explicit form
\begin{align}
\lambda_\phi&=-\frac{b}{\epsilon}\dashint_{\Gamma_0}f^\prime(\phi),
\\
\lambda_u&=-\frac{2\kappa\Lambda \alpha}{R^2}.
\end{align} 

\begin{figure}[ht]
	
					\centering
					\begin{subfigure}{.3\linewidth}
						\centering
						\includegraphics[width=\linewidth]{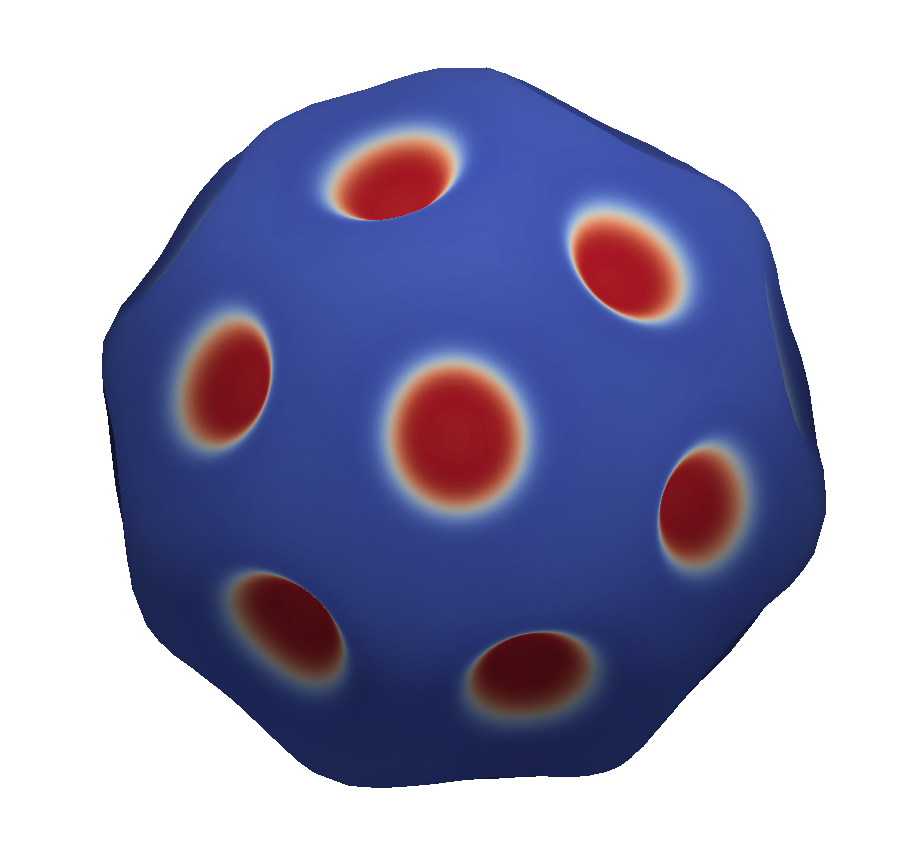}
						\caption{$\Lambda=-10$}
					\end{subfigure}
					\begin{subfigure}{.3\linewidth}
						\centering
						\includegraphics[width=\linewidth]{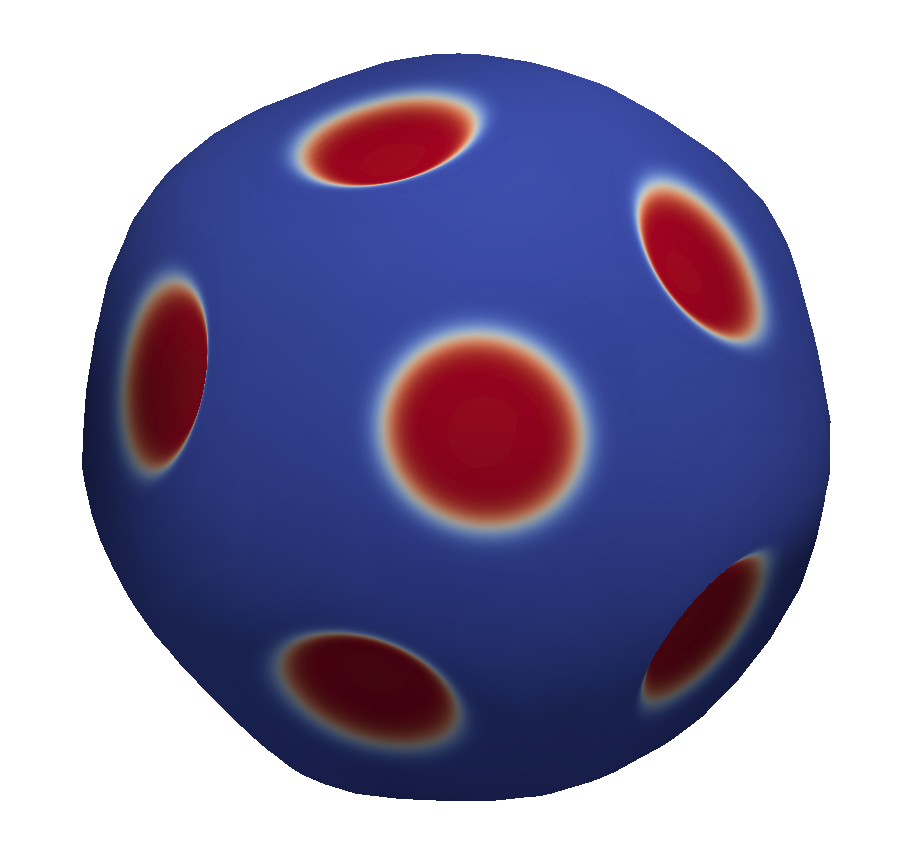}
						\caption{$\Lambda=-5$}
					\end{subfigure}
					\\
					\begin{subfigure}{.3\linewidth}
						\centering
						\includegraphics[width=\linewidth]{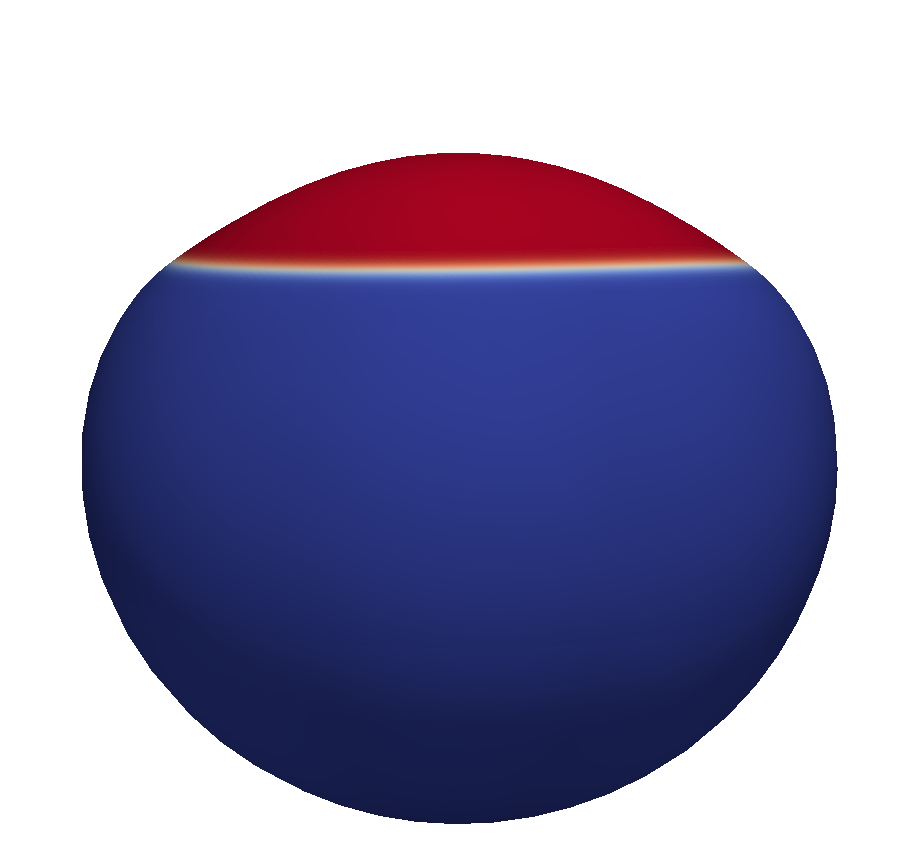}
						\caption{$\Lambda=-1$}
					\end{subfigure}
					\begin{subfigure}{.3\linewidth}
						\centering
						\includegraphics[width=\linewidth]{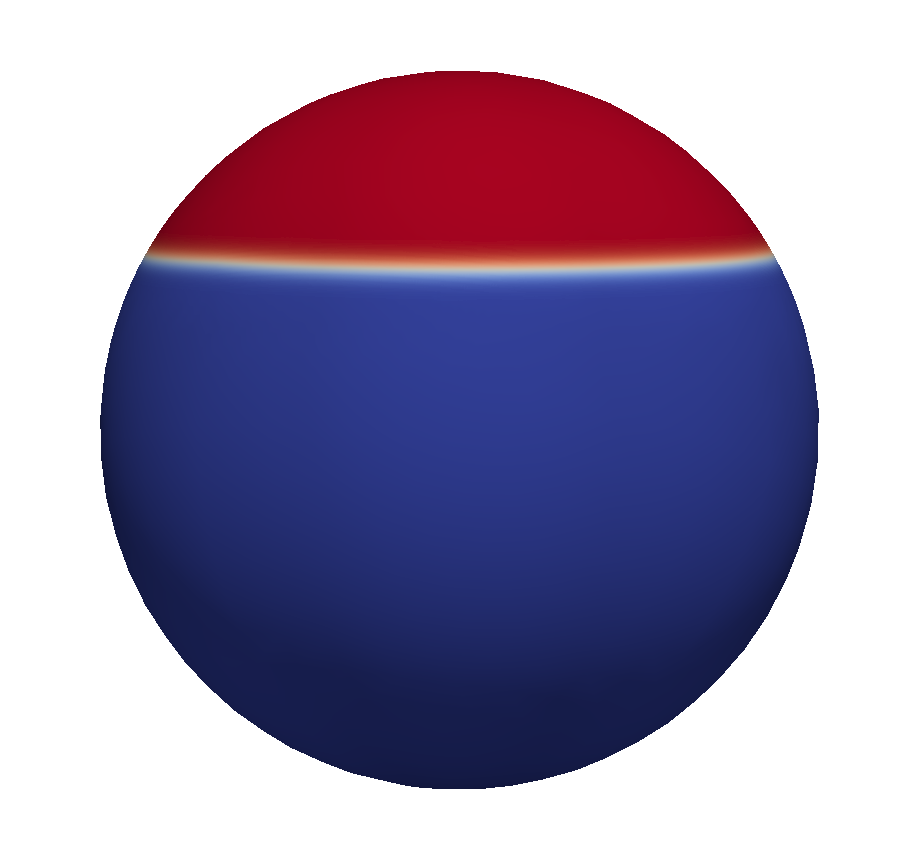}
						\caption{$\Lambda=0$}
					\end{subfigure}
					\begin{subfigure}{.3\linewidth}
						\centering
						\includegraphics[width=\linewidth]{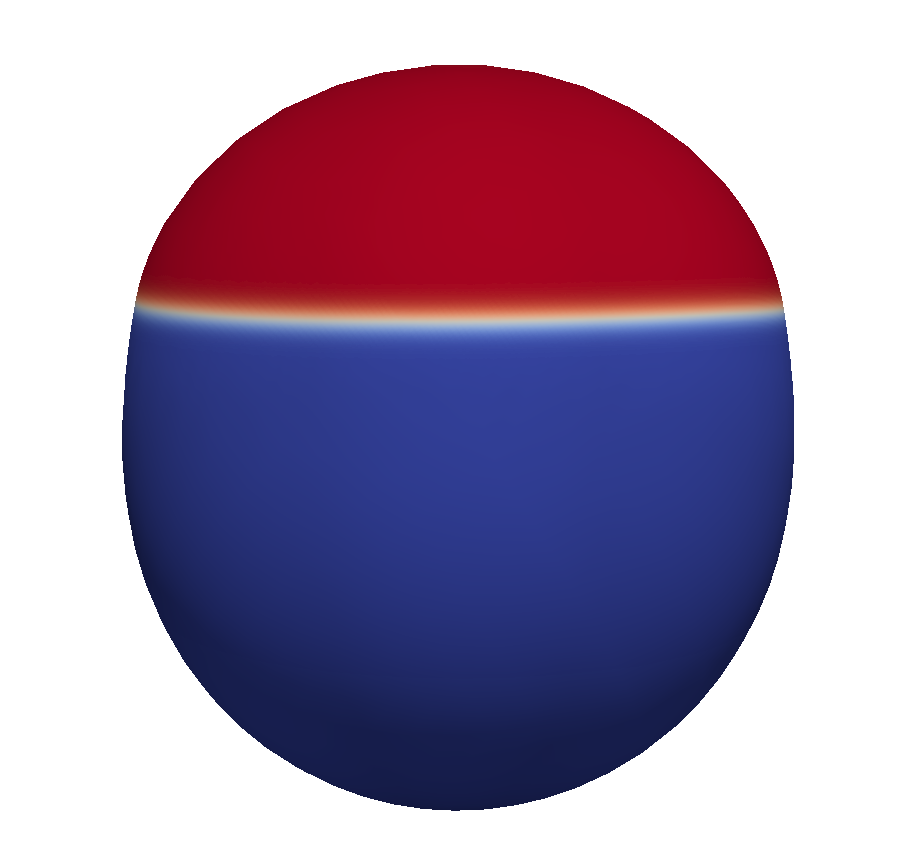}
						\caption{$\Lambda=1$}
					\end{subfigure}
					\begin{subfigure}{.3\linewidth}
						\centering
						\includegraphics[width=\linewidth]{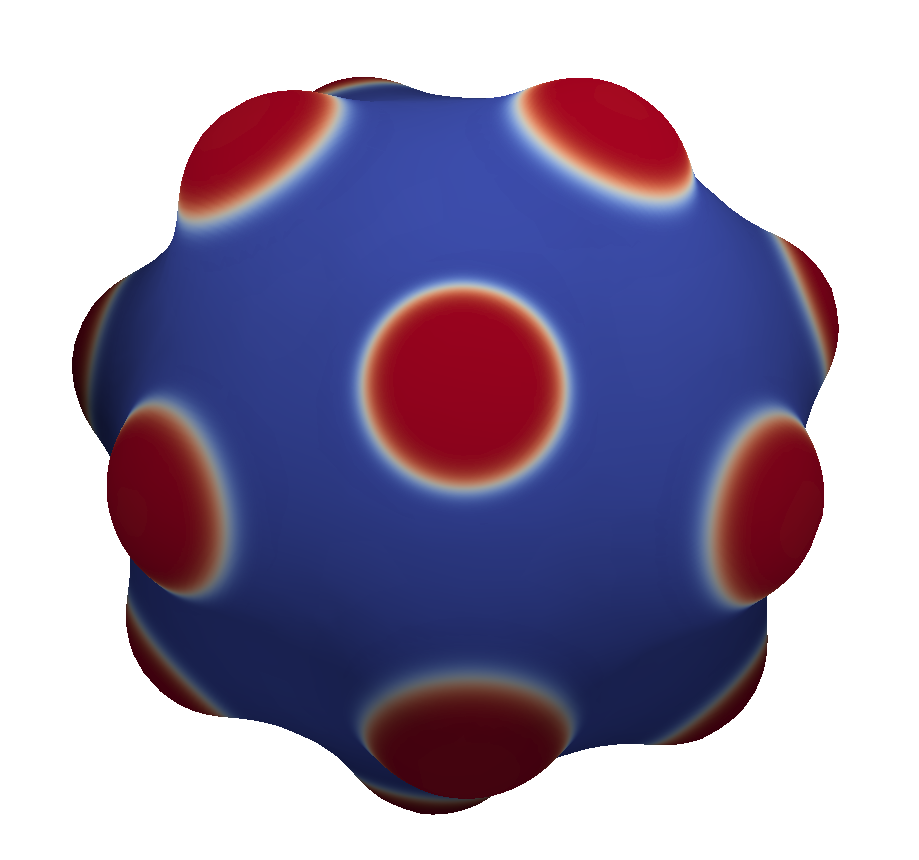}
						\caption{$\Lambda=5$}
					\end{subfigure}
					\begin{subfigure}{.3\linewidth}
						\centering
						\includegraphics[width=\linewidth]{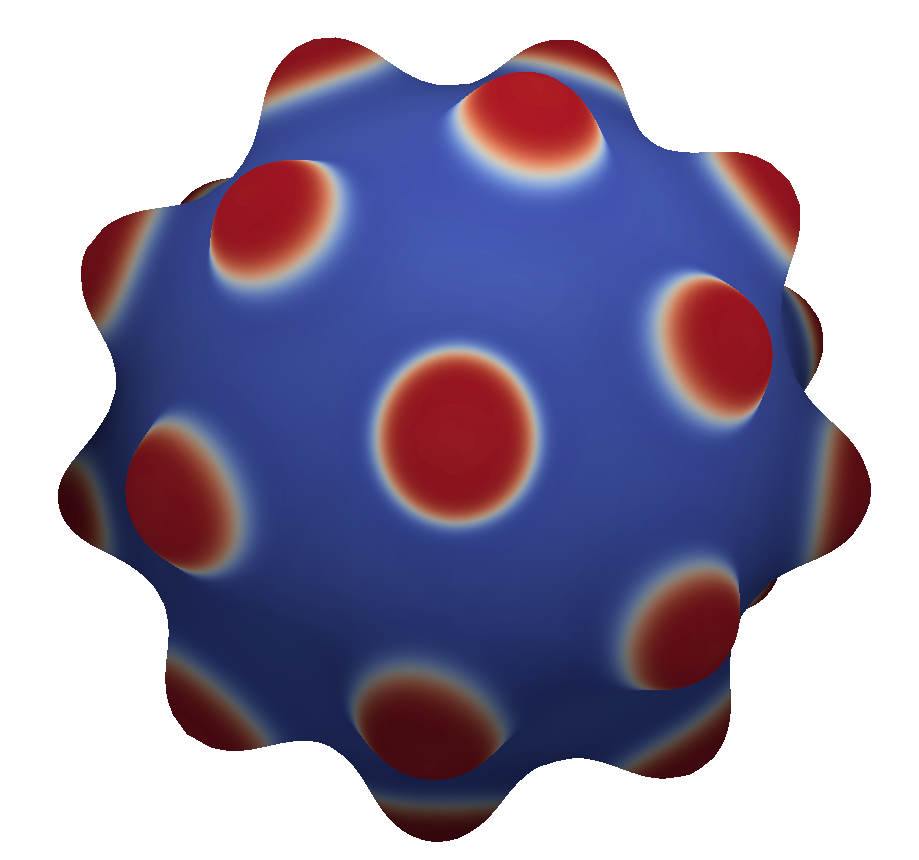}
						\caption{$\Lambda=10$}
					\end{subfigure}
					\caption{Effect of varying the coupling coefficient $\Lambda$. The simulations show the two phases of the order parameter $\phi$. In each case the surface is deformed according to u. For illustration purposes $\rho=1$, whereas in reality it would be much smaller.}
				\label{fig:Coupling}\end{figure}

 DUNE software  with a PYTHON module (c.f. \cite{DedNol18})  was used to implemement a surface FEM method, \cite{DziEll13-a},  for \eqref{AC PDE}. Some examples of almost stationary solutions are given in Figure \ref{fig:Coupling}. Additional details and results for this model will be given in a work in preparation by Elliott, Hatcher and Stinner.
\\
\\

\section{Outlook}\label{Outlook}
 Well posedness and numerical analysis of the problems described in Sections 4 and 5 with detailed proofs are the subjects of works in preparation by the authors. These works also include geometric  motion forced by surface deformation for the sharp interface limit of the  phase field problem and the movement of  point particles on the surface.

\section*{Acknowledgements}
The research of C.M.~Elliott was partially supported by the Royal Society via a Wolfson Research Merit Award.
L. Hatcher and P. Herbert  were supported by the UK Engineering and Physical Sciences Research Council (EPSRC) Grant EP/H023364/1 within the MASDOC Centre for Doctoral Training.

\bibliographystyle{siam}
\bibliography{LibraryBibdeskRefs}
\end{document}